\documentclass{amsart}

\title{A Finite Graph Approach to the Probabilistic Hadwiger-Nelson Problem}
\author{Haydn Gwyn, Jacob Stavrianos}

\usepackage{parskip}
\usepackage[utf8]{inputenc}
\usepackage[letterpaper,portrait,margin=1in]{geometry}
\usepackage{amsmath}
\usepackage{mathtools}
\usepackage{bm}
\usepackage{amssymb}
\usepackage{fancyhdr}
\usepackage{tikz}
\usepackage{hyperref}
\usepackage{mathrsfs}
\usepackage{enumitem}
\usepackage{subfigure}
\usepackage{theorems}
\usepackage{comment}
\usepackage{bbm}

\newcommand{\abs}[1]{\left| #1 \right|}

\newcommand{\C}{\mathbb{C}}
\newcommand{\R}{\mathbb{R}}
\newcommand{\Rnn}{\mathbb{R}_{\geq 0}}
\newcommand{\Q}{\mathbb{Q}}

\newcommand{\Z}{\mathbb{Z}}

\newcommand{\N}{\mathbb{N}}

\newcommand{\E}{\mathbb{E}}

\DeclareMathOperator{\re}{Re}
\DeclareMathOperator{\im}{Im}

\DeclareMathOperator{\pr}{Pr}
\DeclareMathOperator{\cis}{cis}

\newenvironment{hdefn}{\vspace{10pt}\begin{defn}}{\end{defn}\vspace{5pt}}
\newenvironment{hlem}{\vspace{10pt}\begin{lem}}{\end{lem}\vspace{5pt}}
\newenvironment{hthm}{\vspace{10pt}\begin{thm}}{\end{thm}\vspace{5pt}}
\newenvironment{hconj}{\vspace{10pt}\begin{conj}}{\end{conj}\vspace{5pt}}
\newenvironment{hcor}{\vspace{10pt}\begin{cor}}{\end{cor}\vspace{5pt}}

\newtheorem{manualtheoreminner}{Theorem}
\newenvironment{manualtheorem}[1]{
  \renewcommand\themanualtheoreminner{#1}
  \manualtheoreminner
}{\endmanualtheoreminner}

\newtheorem{manualdefninner}{Definition}
\newenvironment{manualdefn}[1]{
  \renewcommand\themanualdefninner{#1}
  \manualdefninner
}{\endmanualdefninner}

\newcommand{\eqd}
{
\begin{tikzpicture}[xscale=0.5,yscale=0.5,
bluenode/.style={circle,draw=blue!60, fill=blue!60, thick, minimum size=5mm},
greennode/.style={circle,draw=green!60, fill=green!60, thick, minimum size=5mm},
rednode/.style={circle,draw=red!60, fill=red!60, thick, minimum size=5mm},
whitenode/.style={circle,draw=black!60,thick,minimum size=5mm},]
\node[bluenode] (blue) at (0,0) { };
\node[rednode] (green) at (4,0) { };
\node[whitenode] (red) at (2,3.464) {?};

\draw[-] (blue) -- (red);
\draw[-] (green) -- (red);
\draw[-] (blue) -- (green);
\end{tikzpicture}
}

\newcommand{\moser}{
\begin{tikzpicture}[xscale=0.5,yscale=0.5,
bluenode/.style={circle,draw=blue!60, fill=blue!60, thick, minimum size=5mm},
greennode/.style={circle,draw=green!60, fill=green!60, thick, minimum size=5mm},
rednode/.style={circle,draw=red!60, fill=red!60, thick, minimum size=5mm},
othernode/.style={circle,draw=black!60, thick, minimum size=5mm},]

\node[othernode] (node1) at (0,0) {?};
\node[bluenode] (node2) at (2.9148542155126758,-2.7392745211657745) {};
\node[greennode] (node3) at (-0.914854215512676,-3.8939750595450255) {};
\node[greennode] (node4) at (0.914854215512676,-3.8939750595450255) {};
\node[rednode] (node5) at (-2.9148542155126758,-2.7392745211657745) {};
\node[bluenode] (node6) at (-2,-6.6332495807108) {};
\node[rednode] (node7) at (2,-6.6332495807108) {};

\draw[-] (node1)--(node2);
\draw[-] (node1)--(node3);
\draw[-] (node1)--(node4);
\draw[-] (node1)--(node5);
\draw[-] (node6)--(node5);
\draw[-] (node6)--(node4);
\draw[-] (node7)--(node2);
\draw[-] (node7)--(node3);
\draw[-] (node4)--(node5);
\draw[-] (node2)--(node3);
\draw[-] (node6)--(node7);
\end{tikzpicture}}

\newcommand{\isbellhex}[3]{
\hexagon{1.5*#1}{0.866*#2}{#3}
}
\newcommand{\isbell}{
\begin{tikzpicture}[xscale=0.6,yscale=0.6]
\isbellhex{0}{0}{black}
\isbellhex{0}{2}{green}
\isbellhex{0}{4}{yellow}
\isbellhex{0}{6}{red}
\isbellhex{1}{1}{red}
\isbellhex{1}{3}{magenta}
\isbellhex{1}{5}{blue}
\isbellhex{1}{7}{cyan}
\isbellhex{2}{0}{blue}
\isbellhex{2}{2}{cyan}
\isbellhex{2}{4}{black}
\isbellhex{2}{6}{green}
\isbellhex{3}{1}{green}
\isbellhex{3}{3}{yellow}
\isbellhex{3}{5}{red}
\isbellhex{3}{7}{magenta}
\isbellhex{4}{0}{red}
\isbellhex{4}{2}{magenta}
\isbellhex{4}{4}{blue}
\isbellhex{4}{6}{cyan}
\isbellhex{5}{1}{cyan}
\isbellhex{5}{3}{black}
\isbellhex{5}{5}{green}
\isbellhex{5}{7}{yellow}
\end{tikzpicture}}

\newcommand{\colorpatch}{
\begin{tikzpicture}[xscale=0.5,yscale=0.5]
\filldraw[-,color=black!20, fill=black!20] (0,0) rectangle (14,3*1.732+2);
\foreach \j in {1,4.464} {
\foreach \i in {1,3,5,7,9,11,13} {
\filldraw[-,color=black!60, fill=red!60] (\i,\j) circle (0.5);}}
\foreach \j in {2.732,3*1.732+1} {
\foreach \i in {2,4,6,8,10,12} {
\filldraw[-,color=black!60, fill=red!60] (\i,\j) circle (0.5);}}
\filldraw [-,color=black!60, fill=red!60,domain=-90:90] plot ({0.5*cos(\x)}, {0.5*sin(\x)+3*1.732+1});
\filldraw [-,color=black!60, fill=red!60,domain=-90:90] plot ({0.5*cos(\x)}, {0.5*sin(\x)+1*1.732+1});
\filldraw [-,color=black!60, fill=red!60,domain=90:270] plot ({14+0.5*cos(\x)}, {0.5*sin(\x)+3*1.732+1});
\filldraw [-,color=black!60, fill=red!60,domain=90:270] plot ({14+0.5*cos(\x)}, {0.5*sin(\x)+1*1.732+1});
\foreach \i in {0,2,4,6,8} {
\draw[-] (\i+0.4226,0)--(\i+4.5773,3*1.732+2);}
\draw[-] (10.4226,0)--(14,6.1962);
\draw[-] (12.4226,0)--(14,2.7321);
\draw[-] (0,3*1.732+2-4.464)--(2.5773,3*1.732+2);
\draw[-] (0,3*1.732+1)--(0.5773,3*1.732+2);
\foreach \i in {0,1,2,3} {
\draw[-] (0,1+\i*1.732)--(14,1+\i*1.732);}
\end{tikzpicture}}

\newcommand{\bars}{
\begin{tikzpicture}[xscale=0.4,yscale=0.4]
\filldraw[color=red!60, fill=red!60, thick] (0,0) rectangle (2,10);
\filldraw[color=blue!60, fill=blue!60, thick] (2,0) rectangle (4,10);
\filldraw[color=red!60, fill=red!60, thick] (4,0) rectangle (6,10);
\filldraw[color=blue!60, fill=blue!60, thick] (6,0) rectangle (8,10);
\filldraw[color=red!60, fill=red!60, thick] (8,0) rectangle (10,10);
\filldraw[color=blue!60, fill=blue!60, thick] (10,0) rectangle (12,10);

\draw[color=black,very thick] (1.8453,4)--(4.1547,4)--(3,6)--cycle;
\draw[color=black,very thick] (6,6)--(6,8.3094)--(8,7.1547)--cycle;
\draw[color=black,very thick] (9.701291823869893, 2.5785116521567444)--(7.398348300912007, 2.7510934107514964)--(8.400359875218101, 0.6703949370917599)--cycle;
\end{tikzpicture}}

\newcommand{\hexagon}[3]{
\filldraw[-,color=#3!60,fill=#3!60] (#1+1,#2)--(#1+0.5,#2+0.866)--(#1-0.5,#2+0.866)--(#1-1,#2)--(#1-0.5,#2-0.866)--(#1+0.5,#2-0.866)--cycle;
}
\newcommand{\hexthree}{
\begin{tikzpicture}[xscale=0.5,yscale=0.5]
\foreach \x in {0,1,2} {\hexagon{3*\x}{0}{blue}}
\foreach \x in {0,1,2} {\hexagon{3*\x+1.5}{0.866}{red}}
\foreach \x in {0,1,2} {\hexagon{3*\x}{1.732}{green}}
\foreach \x in {0,1,2} {\hexagon{3*\x+1.5}{2.598}{blue}}
\foreach \x in {0,1,2} {\hexagon{3*\x}{3.464}{red}}
\foreach \x in {0,1,2} {\hexagon{3*\x+1.5}{4.33}{green}}
\foreach \x in {0,1,2} {\hexagon{3*\x}{5.196}{blue}}
\foreach \x in {0,1,2} {\hexagon{3*\x+1.5}{6.062}{red}}
\end{tikzpicture}}

\newcommand{\hexfour}{
\begin{tikzpicture}[xscale=0.5,yscale=0.5]
\foreach \y in {0,2} {
\foreach \x in {0,1} {\hexagon{3*\y}{1.732*2*\x}{red}}
\foreach \x in {0,1} {\hexagon{3*\y}{1.732*2*\x-1.732}{blue}}
\foreach \x in {0,1} {\hexagon{1.5+3*\y}{1.732*2*\x+0.866}{green}}
\foreach \x in {0,1} {\hexagon{1.5+3*\y}{1.732*2*\x-0.866}{yellow}}}
\foreach \x in {0,1} {\hexagon{3}{1.732*2*\x}{blue}}
\foreach \x in {0,1} {\hexagon{3}{1.732*2*\x-1.732}{red}}
\foreach \x in {0,1} {\hexagon{1.5+3}{1.732*2*\x+0.866}{yellow}}
\foreach \x in {0,1} {\hexagon{1.5+3}{1.732*2*\x-0.866}{green}}
\end{tikzpicture}}

\begin{document}
\begin{abstract}
We advance a probabilistic approach to the Hadwiger-Nelson problem initially developed by the Polymath16 project, in particular relating the approach to finite unit-distance graphs. We define the numerical \textit{badness} of a given $k$-coloring of the plane to be the probability that a randomly chosen unit-distance edge is monochromatic under the coloring, and we provide lower bounds on the badness of arbitrary $k$-colorings using a probabilistic technique relating to finite graphs. The contrapositive of the resulting bounds lets us compute lower bounds on the order of non $k$-colorable unit-distance graphs, improving bounds produced by Pritikin and the Polymath16 project in the $k = 4$ and $k = 5$ cases. Additionally, we make partial progress on a probabilistic analog of the de Bruijn-Erd\H{o}s compactness theorem.
\end{abstract}
\maketitle
\section{Background}
\subsection{Problem Statement}
The \textit{chromatic number} $\chi(G)$ of a graph $G = (V, E)$ is the minimum number of colors $k$ necessary such that each vertex $v \in V$ can be assigned a color in $\{1, 2, \dots k\}$ with the property that any two vertices connected by an edge are assigned different colors.

The \textit{chromatic number of the plane} (CNP) is the minimum number of colors $\chi$ such that the entire plane can be colored with $\chi$ colors in such a way that no two points exactly one unit apart have the same color. \cite{colorbook}. Equivalently, this is the chromatic number of the graph whose vertex set is all points in the plane and whose edge set is the set of pairs of points one unit apart.

The problem of computing the exact value of $\chi$ is referred to as the Hadwiger-Nelson problem, or CNP, and was proposed by Ed Nelson in 1950 \cite{colorbook}. CNP remains an open problem.

\subsection{Classical Results}

A \textit{unit-distance graph} is a graph that can be embedded into the plane with only \textit{unit-distance edges}, or edges connecting vertices distance one apart. Clearly $\chi$ (the chromatic number of the plane) is at least the chromatic number of any finite unit-distance graph. This leads to some simple lower bounds on $\chi$:
\begin{figure}[!ht]
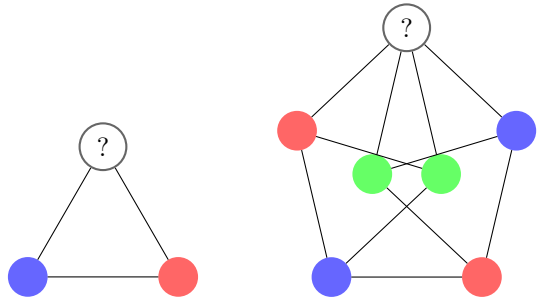

\centering
\begin{tabular}{c c c c}\eqd&~&~&\moser\end{tabular}
\caption{An equilaterial triangle (left) with $\chi(G)=3$ and a Moser Spindle (right) with $\chi(G)=4$.}
\end{figure}

The graphs on the left and right of Figure 1 require three and four colors, respectively, for a valid coloring. This proves that $\chi \geq 4$.

In 1960, Isbell published a hexagonal tiling-based 7-coloring of the plane with each hexagon monochromatic and with diameter slightly less than 1 \cite{colorbook}, thereby showing that $\chi \le 7$:
\begin{figure}[!ht]
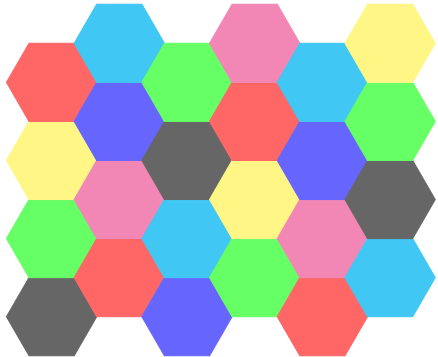

\centering
\isbell
\caption{Isbell's 7-coloring of the plane.}
\end{figure}

The above lower bounds gain a further significance due to a result of de Bruijn and Erd\H{o}s:

\begin{hthm}[de Bruijn, Erd\H{o}s]
The chromatic number of a graph $G$ is equal to the maximum chromatic number of all finite subgraphs $H \subseteq G$.\footnote{This result relies on the Axiom of Choice (AC), as such we prove all results in this paper assuming AC.}
\end{hthm}

In particular, this means that $\chi$ is exactly equal to the maximum chromatic number of all finite unit-distance graphs. As such, if for example $\chi > 4$, then there must exist a non-4-colorable finite unit-distance graph.

Such a graph was found in April 2018 by de Grey, who constructed a 1581-vertex graph with chromatic number five \cite{degrey}. Although this graph was reduced to 553 vertices by Heule, this reduced graph is still quite large, suggesting that non $5$ or $6$-colorable graphs might be extremely large, assuming they exist.

\subsection{Summary of Results}

The present paper expands on the probabilistic approach to CNP developed by the recent Polymath16 Project \cite{probframe}. We work toward a probabilistic analogue of the de Bruijn-Erd\H{o}s theorem.

We define the ``badness'' of a given coloring to be the proportion of unit-distance edges (pairs of points in the plane) that are monochromatic (both endpoints are assigned the same color). We prove that the badness of any $k$-coloring is bounded below by the reciprocal of the number of edges in any non-$k$-colorable unit-distance graph. That is, if a unit-distance graph $G$ is not $k$-colorable and has $E$ edges, then the badness of any $k$-coloring of the plane is at least $\frac1E$. As a corollary, the existence of a $k$-coloring with badness zero (one in which monochromatic edges are sparse) implies that $\chi \leq k$.

As a further corollary of the $\frac1E$ bound, we derive lower bounds on the order of non-$k$-colorable unit-distance graphs from $k$-colorings with low badness. To this end, we use previous research and computer optimization to find colorings with small badness.

Our bounds improve on previous research of Pritikin \cite{pritikin} and the Polymath Project in the four-color and five-color cases. They are summarized in the following table:

\begin{center}\begin{tabular}{c | c | c}
    $k$ &  Lower Bound on $\abs{V}$ & Lower Bound on $\abs{E}$\\
    \hline
    4 & $22$ & $98$\\
    
    5 & $32$ & $178$\\ 
\end{tabular}\end{center}

The bounds on $\abs{E}$ come directly from colorings with low badness. The bounds on $\abs{V}$ come from the bounds on $\abs{E}$ and the known bound
\[\abs E < \abs V ^{3/2}\]
on unit-distance graphs proven by Erd\H{o}s \cite{n32bound}. Naturally, these vertex bounds are weaker than the edge bounds.

Additionally, this paper seeks to strengthen our theorem providing lower bounds on badness from finite graphs into an equality. This equality would be a complete probabilistic analog of the de Bruijn-Erdos theorem. We prove partial results on this equality, but the general statement remains a conjecture.

\section{Probabilistic Approach}

\subsection{Preliminaries}
For ease of notation, we will use $\C$ to refer to the ``plane graph'': a graph with a vertex at every complex number $z$, with vertices $z_1$ and $z_2$ connected by an edge if and only if $\vert z_1 - z_2 \vert = 1$. Further, we define $E(G)$ to be the set of edges in a graph $G$.

In order to rigorize the notion of badness for an arbitrary coloring, we construct a probability measure over the group $E(2)$ of Euclidean isometries\footnote{Not to be confused with the notation $E(G)$ for the edge set of a graph.}. This requires the following definitions:

A \textit{Euclidean isometry} of the plane is a transformation on the plane that preserves the Euclidean distance between any pair of points.

An \textit{amenable} group $G$ is a group for which there exists a finitely additive\footnote{Note that the measure is not necessarily countably additive, which is why we do not refer to it as a probability measure.} measure $\mu$ on $G$ that is invariant under left group action and satisfies $\mu(G)=1$. That is, for any $g \in G$ and $S \subseteq G$, $\mu(S) = \mu(gS)$, and also
\[\mu\Big(\bigcup_{i=1}^n S_i\Big) = \sum_{i=1}^n \mu(S_i)\]
for disjoint $S_1, \cdots, S_n$ \cite{juschenko}. Notably, if $G$ is amenable, then it in fact admits a finitely additive measure $\mu$ invariant on \textit{both sides}, so that $\mu(S)=\mu(gS)=\mu(Sg)$ \cite{lramen}.

\subsection{Formalizing Random Colorings} \label{sec:defpk}
\subsubsection{Notions of Colorings}

Throughout the paper, we discuss multiple different objects under the label of ``coloring''. We define them rigorously here:

\begin{hdefn}
A \textit{k-coloring} with $k \in \N$ colors of a graph $G = (V, E)$ is any function $c: V \to \{1, 2, \dots k\}$. We denote the set of all $k$-colorings (colorings with $k$ colors) of $G$ by $C_k(G)$.
\end{hdefn}
\begin{hdefn}
A \textit{valid} $k$-coloring of $G$ is a coloring such that, for all $v_1, v_2 \in V$ connected by an edge, we have $c(v_1) \neq c(v_2)$.
\end{hdefn}

\begin{hdefn}
A \textit{random} $k$-coloring of $G$ is a random variable $c_r$ over the set of $k$-colorings $C_k(G)$. We require this random variable to be defined everywhere and finitely additive: there exists a function $P(S)$ representing the probability that $c_r \in S$ for all $S \subseteq C_k$, and $P(S_1) + P(S_2) = P(S_1 \cup S_2)$ for disjoint $S_1, S_2$. We denote the set of all random $k$-colorings of $G$ with $C_{r_k}(G)$.
\end{hdefn}

We further define $C_k \equiv C_k(\C)$ and $C_{r_k} \equiv C_{r_k}(\C)$ for ease of notation.

\subsubsection{Randomizing a Fixed Coloring}

Consider a coloring of the plane $c \in C_k(\C)$. We would like to compute the ``badness'' of this coloring by finding the probability that a randomly selected unit-distance pair of points in the plane are assigned the same color by $c$. As such, we need to define a probability distribution over unit-distance pairs of points.

Since $E(2)$ (the group of Euclidean isometries of the plane) is amenable, we can define a finitely additive measure $\mu_T$ on $E(2)$ invariant under group action (Euclidean isometries) on both sides with $\mu(E(2))=1$. Let $T^* \in E(2)$ be a random isometry chosen according to the $\mu_T$ measure (such that $\pr[T^* \in S] = \mu_T(S)$). Then $c^* = c \circ (T^*)^{-1}$ is a random coloring isometric to $c$ for any given value of $T^*$. Additionally, for each $z \in \C$, $c^*(z)$ is a random variable over the set of colors $\{1, 2, \dots k\}$.

We use $(T^*)^{-1}$ so that the random variable $c^*(z)$ will be invariant under isometries of the plane. We demonstrate this fact below:
\begin{align*}
&\phantom{=} \pr[c^*(T(z)) = i]\\
&= \mu_T(\{T^* \mid c \circ (T^*)^{-1} \circ T (z) = i\})\\
&= \mu_T(\{T^* \mid c \circ (T^{-1} \circ T^*)^{-1} (z) = i\}\\
&= \mu_T(\{T \circ X \mid c \circ X^{-1} (z) = i\} \hspace{10pt}\text{(letting $X = T^{-1} \circ T^*$)}\\
&= \mu_T(\{X \mid c \circ X^{-1}(z) = i\})\\
&= \pr[c^*(z) = i]
\end{align*}

With this in mind, we define the ``star operator'':
\begin{hdefn}\label{defn:star}
Consider a $k$-coloring of the plane $c \in C_k$. From $c$, we define a random $k$-coloring $c^* \in C_{r_k}$ by
\[c^* = c \circ (T^*)^{-1}\]
with $T^*$ chosen according to some finitely additive left-invariant measure $\mu_T$ over $E(2)$.\footnote{Interestingly, $E(3)$, the group of Euclidean isometries in three dimensions (i.e.\@ of $\R^3$), is not amenable, which means that the star operator cannot be extended naively to three dimensions.}
\end{hdefn}

As demonstrated above, we have that for any $z \in \C$ and $T \in E(2)$,
\[\pr[c^*(z) = i] = \pr[c^*(T(z)) = i]\]
for each $i \in \{1, \ldots, k\}$.

\subsubsection{Evaluating a Random Coloring}

We now define a metric for the ``badness'' of a $k$-coloring $c$, which intuitively is the proportion of unit-distance pairs in $\C$ that are monochromatic under $c$.

\begin{hdefn}
For a graph $G = (V, E)$ and some $v_1, v_2 \in V$, $c \in C_k(G)$, $c_r \in C_{r_k}(G)$:

\begin{align*}
b_{c}(v_1, v_2) &= \begin{cases} 
      1, \text{ if } c(v_1) = c(v_2)\\
      0, \text{ otherwise} 
   \end{cases}
\end{align*}

\[b_{c_r}(v_1, v_2) = b_{c_0}(v_1, v_2)\]
where $c_0 \in C_k$ in the definition of $b_{c_r}$ is the value of $c_r$. Thus, $b_{c_r}(v_1, v_2)$ is itself a random variable based on the randomness of $c_r$.
\end{hdefn}

In particular, we can apply this definition to colorings of the plane:

\begin{hdefn}
For $c \in C_k(\C)$,
\[p_k(c) = \mathbb{E}_{c^*}[b_{c^*}(0, 1)]\]
where the expected value is taken over all possible values of $c^*$. Rigorously, the expected value represents the integral of $b_{c^*}(0, 1)$ over $E(2)$, as defined in \hyperref[appendix]{the appendix}.
\end{hdefn}

\noindent $p_k(c)$ is essentially the proportion of unit-distance pairs in the plane that are monochromatic under $c$, or equivalently the probability that an edge chosen randomly 
via the randomness of $T^*$ is monochromatic.

There are several things to note here:

\begin{itemize}
    \item $p_k(c)$ is well-defined for all $c \in C_k$.
    \item $p_k(c) \in [0, 1]$.
    \item $\vert z_1 - z_2 \vert = 1 \implies p_k(c) = \E_{c^*} [b_{c^*}(z_1, z_2)]$, by isometry invariance of $c^*$.
    \item $p_k(c) = \mu_T(\{T^* \mid b_{c \circ (T^*)^{-1}}(0, 1) = 1\})$, thus $p_k(c)$ can be represented directly in terms of the measure over $E(2)$ used to define $c^*$.
\end{itemize}

Now we can compute $p_k(c)$ by taking the measure of the set $\{T \in E(2) \mid b_c(T(0), T(1)) = 1\}$ since the expected value is taken over $E(2)$.

Since the edge $(T(0), T(1))$ becomes any given unit-distance edge for exactly two values of $T$, we can intuitively claim that $p_k(\C)$ measures the proportion of monochromatic unit-distance pairs under the coloring $c$.

We can now define $p_k$, which represents the ``best possible'' $k$-coloring in terms of its $p_k(c)$ value:

\begin{hdefn}
For any number of colors $k$,
\[p_k = \inf_{c \in C_k} p_k(c)\]
\end{hdefn}

\subsection{Relating \texorpdfstring{$p_k$}{pk} to Finite Graphs}

The following lemma will be necessary when considering variables with finitely-additive probability spaces:
\begin{hlem}\label{lem:supmean}

Consider some finitely additive probability space $(X, \Sigma, \mu)$. Let $f: X \to \R$ be a measurable function  (where the reals are equipped with some $\sigma$-algebra, say the Borel $\sigma$-algebra). In other words, $f$ is a real-valued random variable over $X$. Then with
\[\E[f] = \int_X f(x) \, d\mu\]
we have $\inf f \leq \E[f] \leq \sup f$.
\end{hlem}

\begin{proof}
First, we note that because the finitely-additive integral is a linear operator (see \hyperref[appendix]{the appendix} for a rigorous proof), the expected value operator is also linear. We note that
\[\sup f = -\inf(-f)\]
so
\[\E[f] \leq \sup (f) \Leftrightarrow -\sup (f) \leq -\E[f] \Leftrightarrow \inf(-f) \leq \E[-f]\]
thus, since $f$ is measurable if and only if $-f$ is measurable, it suffices to prove that $\inf f \leq \E[f]$ for all random variables $f$.

First, if $\inf f = -\infty$, then the result clearly holds, so suppose $\inf f = a \in \R$. Consider the function $g = f - a$, which is also measurable (since the Borel $\sigma$-algebra is translation-invariant). Since $f \geq a$, we have $g = f-a \geq a-a=0$, so $g$ is non-negative. This means that the integral of $g$, being that it is equal to the supremum of the integrals of the simple functions bounded by it, is non-negative as well. So we have
\[0 \leq \int_X g \, d\mu = \int_X f-a \, d\mu = \int f \, d\mu - \int_X a \, d\mu = \E[f] - a \mu(X)\]
whence $\E[f] \geq a\mu(X) = a$, as desired.
\end{proof}

\subsubsection{Finite Graph Bound on $p_k$}
We first offer a definition for convenience of notation.
\begin{hdefn}
For a given coloring $c \in C_k(G)$ and edge $e \in E(G)$, we define $b_c(e) = b_c(v_1, v_2)$, where $v_1$ and $v_2$ are the vertices in $G$ connected by $e$. We similarly define $b_{c_r}(e)$ for random colorings $c_r \in C_{r_k}(G)$.
\end{hdefn}

Now, we develop our first theorem relating $p_k$ to finite graphs:

\begin{hlem}\label{lem:lower}
If a given finite unit-distance graph $G = (V, E)$ is not $k$-colorable, then $p_k \geq \frac{1}{\abs{E(G)}}$.
\end{hlem}

\begin{proof}
Consider an arbitrary embedding (in which all edges have length 1) of $G$ into $\C$. Now, consider an arbitrary coloring $c \in C_k$. We define the random variable
\[B = \sum_{e \in E(G)} b_{c^*}(e)\]

By linearity of expectation and isometry invariance, we have $\mathbb{E}[B] = \abs{E(G)}*p_k(c)$.
A case where $B = 0$ implies a valid $k$-coloring of $G$, so we must have $B \geq 1$ for all values of $c^*$. Thus, $\mathbb{E}[B] \geq 1$ by Lemma \ref{lem:supmean}.

Combining equations, we get $p_k(c) \geq \frac {1} {\abs{E(G)}}$ for all $c \in C_k$. Taking an infimum over colorings $c$, we obtain $p_k \geq \frac {1} {\abs{E(G)}}$.
\end{proof}

\subsubsection{Relating $p_k$ to $\chi(\C)$}

Using our new finite-graph machinery, we specifically analyze the case where $p_k = 0$:

\begin{hthm}\label{thm:lower}
If $\chi(\C) > k$, then $p_k \geq \frac {1} {N}$ for some $N \in \Z^+$.
\end{hthm}
\begin{proof}
Assume that $\chi(\C) > k$, equivalently that the real plane is not $k$-colorable. By the de Bruijn-Erd\H{o}s theorem, there exists a finite unit-distance graph $G = (V, E)$ such that $\chi(G) > k$. Now, apply Lemma \ref{lem:lower} and set $N = \abs{E(G)}$.
\end{proof}
\begin{hcor}\label{cor:zero}
If $p_k = 0$, then $\chi(\C) \leq k$.
\end{hcor}

This follows directly from Theorem \ref{thm:lower}. Hence, if we find a $k$-coloring with $p_k(c) = 0$, then we conclude that $\chi(\C) \leq k$ by Corollary \ref{cor:zero}.
\subsection{Relating Finite Graphs to \texorpdfstring{$p_k$}{pk}}

As motivation, we present the following corollary of Lemma \ref{lem:lower}:

\begin{hlem}\label{lem:eps}
If $0 < p_k \le \varepsilon$, then for any finite unit-distance graph $G = (V,E)$ that is not $k$-colorable, we have $\abs{E(G)} \geq \frac{1}{\varepsilon}$.
\end{hlem}

\begin{proof}
Suppose there exists a unit-distance graph $G = (V,E)$ that is not $k$-colorable and has $\abs{E(G)} < \frac {1} {\varepsilon}$. Then, by Lemma \ref{lem:lower}, $p_k \ge \frac {1} {\abs{E(G)}} > \varepsilon$, a contradiction.
\end{proof}

Thus, we can develop useful lower bounds on the size of non-$k$-colorable unit distance graphs by proving upper bounds on $p_k$. We do this by manually finding and evaluating colorings with low $p_k(c)$ value.

\subsubsection{Calculating $p_k$ on well-behaved colorings}

We now develop machinery for computing $p_k(c)$ for a broad class of ``nice'', periodic colorings $c$. As in the definition of $c^*$, we let $\mu_T$ be the measure over $E(2)$ with the property that $\pr[T^* \in S] = \mu_T(S)$. In particular,
\[p_k(c) = \E_{c^*}[b_{c^*}(0, 1)] = \E_{T^* \in E(2)}[c(T^*(0)) = c(T^*(1))] = \mu_T(\{T^* \in E(2) \mid c(T^*(0)) = c(T^*(1))\})\]

For the remainder of this section, let $z_1, z_2 \in \C \setminus 0$ be two fixed, linearly independent complex numbers.

First, we present some definitions for working with periodic colorings:

\begin{hdefn}
We define the \textit{translation group} $L$ of $z_1, z_2$ to be the set of translations of the plane sending $0$ to $mz_1 + nz_2$ for all $m, n \in \Z$, with the group operation of function composition. Note that $L$ is isomorphic to the lattice group generated by $z_1, z_2$.
\end{hdefn}

\begin{hdefn}
We define $R(z_1, z_2) = \{az_1 + bz_2 \mid 0 \leq a, b < 1\}$, the parallelogram with vertices at $0$, $z_1$, $z_2$, and $z_1 + z_2$.
\end{hdefn}

We can now define periodic colorings:
\begin{hdefn}
A \textit{periodic coloring} is a coloring $c \in \C_k$ invariant under group action by the translation group. Equivalently, $c$ has the property that $c(z) = c(z + z_1) = c(z + z_2)$ for all $z \in \C$.
\end{hdefn}

Restricting our focus to periodic colorings, we now only need to analyze a single parallelogram $R(z_1, z_2)$ to entirely represent the coloring. We call $R(z_1, z_2)$ the \textit{principal rectangle}, and we define the set of isometries $P = \{T \mid T(0) \in R(z_1, z_2)\}$ to be \textit{principal isometries}. Note that every $T \in E(2)$ can be uniquely represented as $l \circ p$ for some $ l \in L, p \in P$.

To analyze the principal rectangle, we define a measure over principal isometries $\mu_P: \mathcal{P}(P) \to [0, 1]$ by $\mu_P(S) = \mu_T(L \circ S)$ (where $L \circ S := \{l \circ s \vert l \in L, s \in S\}$). This measure $\mu_P$ satisfies $\mu_P(\emptyset) = 0$, $\mu_P(P) = 1$, and is invariant under left translation but not necessarily arbitrary isometries.

We now prove an equivalence between $p_k$ and $\mu_P$ justifying our definition of $\mu_P$:

\begin{hthm} \label{thm:pkmup}
Let $c$ be a periodic $k$-coloring, and let $S \in P$ be the set of principal isometries mapping the edge $(0, 1)$ to a monochromatic edge under $c$. Then $p_k(c) = \mu_P(S)$.
\end{hthm}

\begin{proof}
First, note that the set of \textit{all} isometries mapping the edge $(0, 1)$ to a monochromatic edge under $c$ is precisely $L \circ S$. This is a consequence of the invariance of $c$ under left group action by $L$.

Thus, by definition of $p_k(c)$ we have
\[p_k(c) = \mu_T(\{\text{all isometries mapping (0, 1) to a monochromatic edge}\}) = \mu_T(L \circ S) = \mu_P(S)\]

which completes the proof.
\end{proof}

Our goal is to show that $\mu_P$ partially corresponds to Jordan measure over $[0, 1] \times [0, 1] \times [0, 2\pi]$. To do so, we define ``intervals'' $I \subseteq P$ and prove the value of $\mu_P(I)$.

Specifically, for $a, b \in [0, 1]$, $\theta_1 \leq \theta_2 \in [0, 2\pi]$, we define $I(a, b, \theta_1, \theta_2)$ to be the set of isometries $\{T \mid T(0) \in R(a*z_1, b*z_2), T(1) - T(0) = \cis(\theta), \theta \in [\theta_1, \theta_2]\}$. Our goal is to prove that $\mu_P(I(a, b, \theta_1, \theta_2)) = ab * \frac {\theta_2 - \theta_1} {2\pi}$, corresponding to a Jordan measure assigning $\mu_P(P) = 1$.

To achieve our goal, we show that $\mu_P(I(a, b, \theta_1, \theta_2))$ scales linearly with each of $a$, $b$, and $\Delta \theta := \theta_2 - \theta_1$.

We begin by showing linearity in $a$ and $b$, which amounts to first demonstrating the result in the case of scaling by an integer by using a geometric argument. This extends to the case of $a,b \in \Q$ applying the integer scaling argument twice (once for numerators and once for denominators), and we then rely on the density of the rationals in $\R$ to show the result for all $a,b \in \R$.

\begin{hlem}\label{lem:abscale}
\[\mu_P(I(a, b, \theta_1, \theta_2)) = ab * \mu_P(1, 1, \theta_1, \theta_2)\]
\end{hlem}

\begin{proof}

First, we show the integer scaling case:
\[\forall n \in \Z^+ \hspace{0.2cm} \mu_P(I(a, b, \theta_1, \theta_2)) = n * \mu_P\left(I\left(\frac {a}{n}, b, \theta_1, \theta_2\right)\right)\]
This is seen by noting that the set $I(a, b, \theta_1, \theta_2)$ is the disjoint union of $n$ translated copies of $I(\frac {a} {n}, b, \theta_1, \theta_2)$. Visually, this is equivalent to lining up $n$ copies of a parallelogram into one longer parallelogram, with one side scaled by a factor of $n$.

Next, we show the rational scaling case by applying the integer case twice:
\[\forall \, \frac pq \in \Q, 0 < \frac{p}{q} \leq 1, \hspace{0.2cm} \mu_P\left(I\left(\frac {p} {q}, b, \theta_1, \theta_2\right)\right) = p \cdot \mu_P \left(I\left(\frac1q,b,\theta_1,\theta_2\right)\right) =
\left(\frac {p}{q}\right)*\mu_P(I(1, b, \theta_1, \theta_2))\]

Note additionally that $\frac pq = 0 \rightarrow \mu_P(I(\frac pq, b, \theta_1, \theta_2)) = 0$, thus the claim holds for all $\frac pq \in [0, 1]$.

We can then extend to the real number scaling case:
\[\forall a \in [0, 1], \hspace{0.2cm} \mu_P(I(a, b, \theta_1, \theta_2)) = a * \mu_P(I(1, b, \theta_1, \theta_2))\]

This claim follows from noting that the function $\mu_P(I(a,b,\theta_1,\theta_2))$ is strictly increasing in $a$ and observing that the claim holds when $a$ is rational. Because the rationals are dense in the reals, the result holds.

Applying the same result with $a$ and $b$ switched, we get:
\[\forall a, b \in [0, 1], \hspace{0.2cm} \mu_P(I(a, b, \theta_1, \theta_2)) = a * \mu_P(I(1, b, \theta_1, \theta_2)) = ab * \mu_P(I(1, 1, \theta_1, \theta_2))\]

which completes the proof.
\end{proof}

We now show linearity in $\Delta \theta$ given $a=b=1$. In particular, we show that $\Delta \theta$ determines $\mu_P$ (that is, shifting $\theta_1$ and $\theta_2$ by a constant leaves $\mu_P$ unchanged), and we then employ a similar argument to Lemma \ref{lem:abscale} to $\Delta \theta$ to complete the claim.

\begin{hlem}\label{lem:thetascale}
\[\mu_P(I(1, 1, \theta_1, \theta_2)) = \frac {\theta_2 - \theta_1} {2\pi}\]
\end{hlem}

\begin{proof}
We show invariance under constant shift in $\theta_1$ and $\theta_2$ by writing out the definitions of $\mu_P$ and $I$. Observe that
\begin{align*}
\mu_P(I(1, 1, \theta_1, \theta_2)) &= \mu_T(L \circ I(1, 1, \theta_1, \theta_2)) \\&= \mu_T(\{T \mid T(1) - T(0)=\cis(\theta), \theta \in [\theta_1, \theta_2]\})
\end{align*}
and
\begin{align*}
\mu_P(I(1, 1, \theta_1 + \theta_s, \theta_2 + \theta_s)) &= \mu_T(L \circ I(1, 1, \theta_1 + \theta_s, \theta_2 + \theta_s)) \\&= \mu_T(\{T \mid T(1) - T(0) = \cis(\theta), \theta \in [\theta_1 + \theta_s, \theta_2 + \theta_s]\})
\end{align*}

Since $\mu_T$ is invariant under arbitrary isometries, we can apply a rotation by $\theta_s$ to the first set to transform it exactly into the second set. Thus, the two sets have the same measure.

The claim $\mu_P(I(1, 1, \theta_1, \theta_2)) = \frac {\Delta \theta} {2\pi}$ can be shown analogously to Lemma \ref{lem:abscale} applied to $\Delta \theta$, specifically through the following claims: 
\begin{enumerate}
    \item $\mu_P(I(1, 1, 0, \theta)) = n*\mu_P(I(1, 1, 0, \frac {\theta} {n}))$ by invariance under angle translation
    \item $\mu_P(I(1, 1, 0, \frac {p} {q} 2\pi)) = \frac {p} {q}*\mu_P(I(1, 1, 0, 2\pi)) = \frac {p} {q}$ by applying the above twice and $\mu_P(I(1, 1, 0, 2\pi)) = 1$
    \item $\mu_P(I(1, 1, 0, \theta)) = \frac {\theta} {2\pi}$ since $\mu_P(I(1, 1, 0, \theta))$ is strictly increasing in $\theta$
    \item $\mu_P(I(1, 1, \theta_1, \theta_2)) = \frac {\theta_2 - \theta_1} {2\pi}$ by angle translation by $-\theta_1$
\end{enumerate}
which completes the proof.
\end{proof}

Combining our two scaling lemmas, we can now prove the final result:
\begin{hthm}\label{thm:intvlmes}
$\mu_P(I(a, b, \theta_1, \theta_2)) = ab * \frac {\theta_2 - \theta_1} {2\pi}$
\end{hthm}

\begin{proof}
Combine Lemmas \ref{lem:abscale} and \ref{lem:thetascale}.
\end{proof}

Using ideas from integration theory, we now show that $\mu_P$ corresponds to Jordan measure over a much broader collection of sets.

Noting the bijection $\varphi:T \mapsto (\re T(0),\im T(0), \arg(T(1)-T(0)))$, we can parameterize $P$ with the space $P' = [0, 1] \times [0, 1] \times [0, 2\pi]$. We can consider the standard Jordan measure\footnote{Technically ``Jordan measure'' is not a measure, since its underlying algebra is only finitely additive} $\mu_J$ scaled so that $\mu_J(P')=1$. With this established, we state the following result.

\begin{hthm}\label{thm:mupeval}
For all $S \subseteq P$ such that $\varphi(S)$ is Jordan measurable, we have
\[\mu_P(S) = \mu_J(\varphi(S))\]
\end{hthm}

\begin{proof}
Note that ``rectangles'' in the Jordan measure sense correspond to intervals (translated copies of $I(a, b, \theta_1, \theta_2)$). Letting $S$ be an interval, the claim is proven by Theorem \ref{thm:intvlmes}.

Thus, if $\varphi S$ is Jordan-measurable, then $\mu_P(S) = \mu_J(\varphi S)$, which proves the claim.
\end{proof}

With this result established, we freely write $\mu_J(S)$ in place of $\mu_J(\varphi(S))$, recalling the equivalence. Now, with our machinery in place, we can finally make statements about $p_k(c)$:

\begin{hthm}\label{thm:pkceval}
Let $c$ be a periodic coloring. Then $p_k(c) = \int_P b_c(e) d\mu_J(e)$ if the integral exists.
\end{hthm}

\begin{proof}

By Theorem \ref{thm:pkmup}, we have $p_k(c) = \mu_P(S)$, where $S \subseteq P$ is the set of principal isometries mapping $e_0$ to a monochromatic edge. By Theorem \ref{thm:mupeval}, we have $p_k(c) = \mu_J(S)$. Thus, it remains to prove that the integral on the right-hand side is equal to $\mu_J(S)$.

Note that any integral with respect to Jordan measure is a Riemann integral, and the Riemann integral of the indicator function of $S$ is equal to the Jordan measure of $S$ if either exist \cite{riemjord}, which completes the proof.
\end{proof}

Theorem \ref{thm:pkceval} allows us to compute the value of $p_k(c)$ for a periodic coloring $c$ by taking a Riemann integral over a single period of the coloring. Thus, our goal will be to find periodic colorings that minimize the value of this integral, which we can both computationally approximate and compute exactly.

\subsection{Upper Bounds on \texorpdfstring{$p_k$}{pk}}

To find colorings minimizing the value of $p_k(c)$, we started with well-known colorings of the plane and scaled them by some scaling parameter. We then computationally optimized the scaling parameter value via Monte Carlo approximation of $p_k(c)$.

We will also make use of the following lemma:
\begin{hlem}\label{lem:precur}
$p_{k + 1} \leq (1 - \frac {\pi} {4\sqrt{3}})p_k$
\end{hlem}
\begin{proof}
Consider some $k$-coloring $c$ with $p_k(c) = x$. We will overlay a set $S$ of color $k + 1$ (occupying a $\frac{\pi}{8\sqrt{3}}$ fraction of the plane and adding no monochromatic edges) onto the coloring $c$.

Specifically, we consider the following infinite triangular grid of circles:
\begin{figure}[!ht]
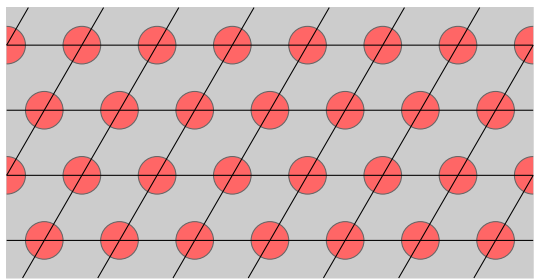

    \centering
    \colorpatch
    \caption{A ``$k+1$-color patch'' (the red circles) to overlay on a $k$-coloring $c$ (the gray background).}
\end{figure}

Each circle has diameter 1, and the distance between 2 of the closest circles is 1, thus no unit-distance edges have both endpoints in $S$.

Note that $S$ is periodic in $z_1 = 2, z_2 = 2e^{i\pi/3}$. With this in mind, we will apply our results on periodic colorings.

We will now define a random $k+1$-coloring. Let $T_P \in P$ be a random variable over isometries such that $T_P$ is always a translation and $T_P(0)$ is uniformly distributed over $R(z_1, z_2)$.

Now, consider the random $k+1$-coloring $c_r$ based on $c$ the random variable $T_P$:
\begin{align*}
c_r(z) &= \begin{cases} 
      k + 1, \text{ if } z \in T_P(S) \\
      c(z), \text{ if } z \notin T_P(S)
   \end{cases}
\end{align*}

Note that if either of $z_1, z_2 \in T_P(S)$, then $b_{c_r}(z_1, z_2) = 0$.

We will show that $\mathbb{E} [p_{k+1}(c_r)] = (1 - \frac {\pi} {4\sqrt{3}}) p_k(c)$ through the following computation:
\begin{align*}
    \mathbb{E}_{c_r}(p_{k + 1}(c_r)) &= \mathbb{E}_{T_P} \mathbb{E}_{T^*} [b_{c_r}(e_0) ] \\
    &= \mathbb{E}_{T_P} \mathbb{E}_{T^*} [\mathbf{1}[T^*(e_0) \cap T_P(S) = \emptyset \land b_c(T^*(e_0)) = 1]] \\
    &= \mathbb{E}_{T_P} \mathbb{E}_{T^*} [\mathbf{1}[T^*(e_0) \cap T_P(S) = \emptyset] * \mathbf{1}[b_c(T^*(e_0)) = 1]]
\end{align*}
Now, define the following subsets of $E(2)$:
\begin{itemize}
    \item $\tau_S = \{T \mid T(e_0) \cap S = \emptyset \}$
    \item $\tau_c = \{T \mid b_c(T(e_0)) = 1\}$
\end{itemize}

We have $\mu_T(\tau_S) = 1 - 2\mu_P(S) = 1 - \frac {\pi} {4\sqrt{3}}$ by Theorem \ref{thm:pkceval} and $\mu_T(\tau_c) = p_k(c)$ by definition of $p_k$.

Given this notation, we can simplify our calculation further:
\begin{align*}
    \mathbb{E}_{c_r}(p_{k + 1}(c_r)) &= \mathbb{E}_{T_P} \mathbb{E}_{T^*} [\mathbf{1}[T^* \in T_P(\tau_S)] * \mathbf{1}[T^* \in \tau_c]]
\end{align*}
Both expected values are defined using at least finitely-additive and thus linear\footnote{see \hyperref[appendix]{the appendix} for proof} integrals, so linearity of expectation holds. After rearranging and invoking the left-amenability of the measure over $T_P$, we can substitute the $\mathbf{1}[T^* \in T_P(\tau_S)]$ term for the summation
\[\frac1N \sum_{i = 1}^N \mathbf{1}[T^* \in T_i \circ T_P(\tau_S)]\]

Upon choosing the $T_i$ "uniformly" over the support of $T_P$ (letting the underlying choice of point $T_i(0) \in R(z_1, z_2)$ approach uniform as $N \to \infty$), the resulting term approaches the constant function $\frac {\mu_T(\tau_S)} {\mu_T(E(2))} = 1 - 2\mu_P(S)$. Making the substitution, we get:

\begin{align*}
    &= \mathbb{E}_{T_P} \mathbb{E}_{T^*} [(\mu_P(S)) * \mathbf{1}[T^* \in \tau_c]] \\
    &= (1 - 2\mu_P(S)) * \mathbb{E}_{T_P} \mathbb{E}_{T^*} [\mathbf{1}[T^* \in \tau_c]] \\
    &= (1 - 2\mu_P(S)) * \mathbb{E}_{T_P} [(p_k(c))] \\
    &= (1 - 2\mu_P(S)) * (p_k(c)) =  (1 - \frac {\pi} {4\sqrt{3}}) p_k(c)\\
\end{align*}

Thus, we conclude that $\mathbb{E}_{c_r}(p_{k + 1}(c_r)) = (1 - \frac {\pi} {4\sqrt{3}}) p_k(c)$. Invoking Lemma \ref{lem:supmean}, there exists some value $c_0$ of $c_r$ yielding $p_{k + 1} \leq p_{k + 1}(c_0) \leq (1 - \frac {\pi} {4\sqrt{3}}) p_k(c)$. This completes the proof.
\end{proof}
\subsubsection{Values for given k}

We now provide periodic $k$-colorings that, combined with Theorem \ref{thm:pkceval}, yield upper bounds on $p_k$ for each of $k \in \{2, 3, 4, 5\}$. For completeness, we have $p_1 = 1$ and $p_k = 0$ for $k \geq 7$.

\begin{itemize}
    \item [$k = 2$:] Let $c$ be a 2-coloring of the plane consisting of alternating stripes of width $\frac {\sqrt{3}} {2}$, each containing all points on their left border.
    
    \begin{figure}[!ht]
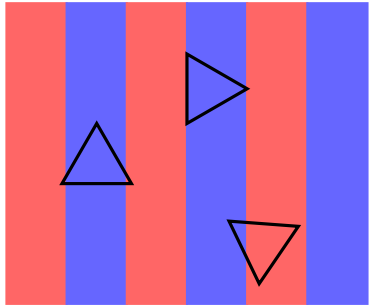

    \centering
    \bars
    \caption{A 2-coloring composed of stripes of width $\frac{\sqrt3}2$ with a $p(c)$-value of $\frac13$.}
    \end{figure}
    
    We can compute $p_2(c) = \frac {1} {3}$. In fact, by applying Lemma \ref{lem:lower} with $G$ an equilateral triangle graph, we find that $p_2 = \frac {1} {3}$ precisely.
    
    \item [$k = 3$:] Let $c$ be the following hexagonal $3$-coloring:
    
    \begin{figure}[!ht]
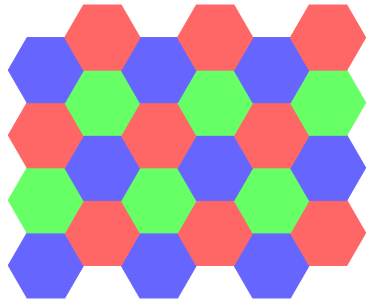

    \centering
    \hexthree
    \caption{A 3-coloring made from a tiling of the plane with hexagons of diameter 1.22.}
    \end{figure}
    
    We can compute $p_3(c) \approx 0.121$.
    
    \item [$k = 4$:] Let $c$ be the following hexagonal $4$-coloring:
    \begin{figure}[!ht]
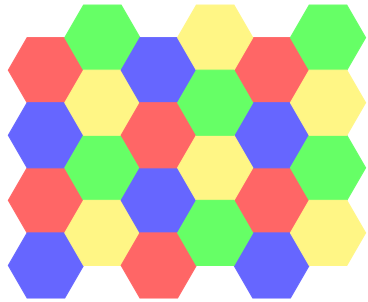

    \centering
    \hexfour
    \caption{A 4-coloring of made from a tiling of the plane with hexagons of diameter 1.13.}
    \end{figure}
    
    We can compute $p_4(c) \approx 0.0102$.
    
    \item [$k = 5$:] To provide an upper bound for $p_5$, we use the upper bound for $p_4$ and Lemma \ref{lem:precur}, yielding the following bound:
    \[p_5 \leq (1 - \frac {\pi} {8\sqrt3})p_4 \leq (1 - \frac {\pi} {8\sqrt3})(0.0101527) \approx .005\]

\end{itemize}

\subsection{Lower Bounds on \texorpdfstring{$k$}{k}-chromatic Graph Size}

Using Lemma \ref{lem:eps}, we can directly convert the upper bounds on $p_k$ into lower bounds on the number of edges in a non-$k$-colorable unit-distance graph.

To derive lower bounds on the number of vertices in such graphs, we use the bound $\vert V \vert > \vert E \vert ^{\frac {2} {3}}$ proven by Erd\H{o}s \cite{n32bound}. This is a relatively weak bound, so the lower bounds on $\vert V \vert$ below are weaker than the lower bounds on $\vert E \vert$.

We summarize our bounds in the table below:
\begin{center}\begin{tabular}{c | c | c | c | c}
    $k$ & Upper Bound on $p_k$ & Lower Bound on $p_k$ & Lower Bound on $\abs{E}$ & Lower Bound on $\abs{V}$ \\
    \hline
    2 & $1/3$ & $1/3$ & 3 & 3 \\
    
    3 & $.121$ & $1/11$ & $9$ & $5$ \\ 
    
    4 & $.0101528$ & $1/2722$ & $98$ & $22$\\
    
    5 & $.00563$ & 0 & $180$ & $32$\\ 

\end{tabular}\end{center}

Notably, the vertex lower bounds of 22 and 32 for $k=4$ and $k=5$ are stronger than the current best known bounds for $k=4$ and $k=5$, which are 13 and 25, respectively.

\section{Finitary Representation of \texorpdfstring{$p_k$}{pk}}
\subsection{Defining a Lower Bound on \texorpdfstring{$p_k$}{pk}}

So far, our only tool for generating lower bounds on $p_k$ is Lemma \ref{lem:lower}. However, this bound is tight only in very specific cases; the lemma can be easily refined by considering multiple monochromatic edges or edge weights. In this section, we present a generalization of the lower-bounding technique to define a stronger lower bound, which we call $q_k$.

\begin{hdefn}
Consider an arbitrary graph $G = (V,E)$ and $k$-coloring $c \in C_k(G)$. Further, consider a nonnegative real-valued weighting function $w: E(G) \to \Rnn$ such that $\sum_{e \in E(G)} w(e)$ is finite.

We now define
\[q_k(G, w) = \inf_c \frac {\sum_{e \in E(G)} w(e)b_c(e)} {\sum_{e \in E(G)} w(e)}\]
\[q_k(G) = \sup_w q_k(G, w)\]
Intuitively, $q_k(G)$ is a weighted average of $b_c$ over $E(G)$, choosing the weights to maximize the amount of ``badness'' necessary in any coloring.

To produce a value comparable to $p_k$, we define
\[q_k = q_k(\C)\]
\end{hdefn}

We will prove that $p_k \geq q_k$ in Section \ref{sec:relpq}, representing a stronger version of Lemma \ref{lem:lower}.

\begin{hlem}
For graphs $G, H$, $G \subseteq H$ implies $q_k(G) \leq q_k(H)$.
\end{hlem}

\begin{proof}
Consider an arbitrary weighting $w_G$ on $E(G)$. We can extend $w_G$ to a weighting $w_H$ on $E(H)$ as follows:

\begin{align*}
w_H(e) &= \begin{cases} 
      w_G(e), \text{ if } e \in E(G)\\
      0, \text{ if } e \notin E(G)
   \end{cases}
\end{align*}
\\
This yields $q_k(G, w_G) = q_k(H, w_H)$, since the weighted averages are the same for all colorings $c$. By taking a supremum over choice of $w_G$, we get $q_k(G) \leq q_k(H)$.
\end{proof}

\subsection{Generalizing \texorpdfstring{$p_k$}{pk} to Arbitrary Graphs} \label{sec:genpk}

To relate $p_k$ and $q_k$ more directly, we now extend our definition of $p_k$, replacing the plane graph $\C$ with an arbitrary graph $G$.

\begin{hdefn}
We define
\[p_k(G, c_r) = \sup_{e \in E(G)} \mathbb{E}[b_{c_r}(e)]\]
where $c_r \in C_{r_k}(G)$ is any random variable $k$-coloring of the graph $G$. And as before, we define

\[p_k(G) = \inf_{c_r} p_k(G, c_r)\]

\end{hdefn}

The equivalence of this definition with our earlier definition of $p_k$ is unclear, so we prove it explicitly:

\begin{hthm}\label{thm:pequiv}
$p_k(\C)$ (as defined above) equals $p_k$ (as defined in Section \ref{sec:defpk}). Thus, the new definition of $p_k$ is an extension of the previous one.
\end{hthm}

\begin{proof}

We have $p_k(\C) \leq p_k$ since $c^*$ is a random coloring for all colorings $c$.
Now, suppose we have some $c_r$ such that $p_k(\C, c_r) = x$. We define a new coloring $c_m = c_r \circ T^* \in C_{r_k}$ as in Definition \ref{defn:star}.

Let $e_0$ be the edge with endpoints 0 and 1. For any fixed $T^* \in E(2)$, we have $\mathbb{E}[b_{c_m}(e_0)] \leq x$. Thus, $\mathbb{E}[b_{c_m}(e_0)] \leq x$, where the expectation is over choice of $T^*$ and $c_r$, by Lemma \ref{lem:supmean}.

Again by Lemma \ref{lem:supmean}, we get that there exists some fixed value $c$ of $c_r$ that yields $\mathbb{E}[b_{c_m}(e_0)] \leq x$. Substituting in for $c_m$, we get $\mathbb{E}[b_c(T^*(e_0))] \leq x$, which is equivalent to $p_k(c) \leq x$.

\end{proof}

\subsection{Relating \texorpdfstring{$p_k$}{pk} and \texorpdfstring{$q_k$}{qk}} \label{sec:relpq}

We can now directly relate $p_k(G)$ and $q_k(G)$ in an analogue of Lemma \ref{lem:lower}:

\begin{hthm}\label{thm:pgeqq}
$p_k(G) \geq q_k(G)$ for all graphs $G$.
\end{hthm}

\begin{proof}
The proof is based on the logic of Lemma \ref{lem:lower}. Assume $p_k(G, c_r) = x$ for some $c_r$. Now, for arbitrary choice of $w$, we define
\[B = \sum_{e \in E(G)} w(e)b_{c_r}(e)\]
From this, we obtain the following bounds:
\[\mathbb{E}_{c_r}[B] \leq \left( \sum_{e \in E(G)} w(e) \right) p_k(G, c_r)\]
\[B \geq \inf_{c_r = c} \sum_{e \in E(G)} w(e)b_c(e)\]
Combining yields
\[p_k(G, c_r) \geq \inf_{c_r = c} \frac {\sum_{e \in E(G)} w(e)b_c(e)} {\sum_{e \in E(G)} w(e)} \geq \inf_c \frac {\sum_{e \in E(G)} w(e)b_c(e)} {\sum_{e \in E(G)} w(e)}\]
Taking a supremum over choice of $w$ yields $q_k(G) \leq x$, from which the theorem follows.
\end{proof}

\begin{hcor}
$p_k \geq q_k$
\end{hcor}

\begin{proof}
Apply Theorem \ref{thm:pgeqq} with $G = \C$ and use Theorem \ref{thm:pequiv} to relate back to the original definition of $p_k$.
\end{proof}

\begin{hthm}\label{thm:finpq}
$p_k(G) = q_k(G)$ for all finite graphs $G$.
\end{hthm}

\begin{proof}
We have $p_k(G) \geq q_k(G)$ by Theorem \ref{thm:pgeqq}, it remains to prove that $p_k(G) \leq q_k(G)$. For the sake of contradiction, assume there exists an $x \in \R$ such that $q_k(G) = x$ and $p_k(G) > x$.

We can represent a given non-variable coloring $c$ as the $|E(G)|$-dimensional vector 
\[\vec{c} = \langle b_c(e_1), b_c(e_2), \dots b_c(e_{\abs{E(G)}})\rangle\]
With this vector representation, we can represent variable colorings as a weighted sum of non-variable colorings:
\[\vec{c_r} = \sum_i m_i \vec{c_i} \hspace{10 pt} \Big(\sum_i m_i = 1\Big)\]
Note that since $G$ is finite, the set $C_k(G)$ is finite, so the above is a finite sum and classifies all possible variable colorings of $G$.

We define $R_1$ as the region of all possible $\vec{c_r}$ and $R_2$ to be the region of the space in which each coordinate is less than or equal to $x$. Both of these regions are convex and, by the $p_k(G) > x$ assumption, are disjoint.

By the Hyperplane Separation Theorem, there must be some hyperplane separating $R_1$ and $R_2$. We write it as $\vec{n} \cdot \vec{v} = y$ for some normal vector $\vec{n}$ and $y \in \R$.

For the hyperplane to be disjoint from $R_2$, all coordinates of $\vec{n}$ must have the same sign or be 0. Without loss of generality let them be nonnegative. Additionally, the vector $\vec{x} = \langle x, x, \dots x \rangle \in R_2$ must be on the $R_2$ side of the plane. Thus, our equation becomes $\vec{n} \cdot (\vec{v} - \vec{x}) = y' > 0$.

Now, we choose the weighting $w(e_i) = n_i$ ($n_i$ denotes the coordinate of $\vec{n}$ corresponding to $e_i$). From the $q_k(G) = x$ assumption and plugging in $w$, we get that there exists some coloring $c$ such that

\[\frac {\sum_{e \in E(G)} w(e)b_c(e)} {\sum_{e \in E(G)} w(e)} \leq x.\]

Rewriting in terms of $\vec{c}$ and substituting $\vec{n}$ for the $w$ terms, we get

\[\vec{n} \cdot \vec{c} \leq \vec{n} \cdot \vec{x}\]

but $\vec{c} \in R_1$ by definition, so this contradicts the hyperplane equation.
\end{proof}

\begin{hconj} \label{conj:pkequalsqk}
$p_k(G) = q_k(G)$ for general graphs $G$. By consequence, $p_k = q_k$.
\end{hconj}

\noindent

A proof of Conjecture \ref{conj:pkequalsqk} would substantially strengthen our results, allowing us to represent $p_k(G)$ directly in terms of finite graphs and to extend our probabilistic methods to higher dimensions. Additionally, the conjecture represents a full probabilistic analogue of the de Bruijn-Erd\H{o}s theorem, establishing an equivalence between ``nicely'' $k$-coloring a graph and ``nicely'' $k$-coloring all its finite subgraphs.
\section*{Appendix: Linearity of the Integral Over Finitely Additive Measures}
\label{appendix}
Here we attempt to clarify some ideas about finitely additive measures. In particular, we demonstrate that linearity of the integral operator continues to hold when we pass into the finitely additive realm. As it happens, our primary use case is functions mapping from $E(2)$, the set of Euclidean isometries of the plane, to the set $\{0,1\} \subset \R$, but for the sake of generality, let $X$ be a set equipped with a $\sigma$-algebra $\Sigma$ and a \textit{finitely} additive measure $\mu$. That is, if $A_1, \cdots, A_n$ are disjoint elements of $\Sigma$, then
\[\mu \left( \bigcup_{i=1}^n A_i \right) = \sum_{i=1}^n \mu(A_i)\]
but if we replace $n$ with $\infty$, the result no longer necessarily holds. We can now define the integral analogously to how it is defined in classical measure theory, using the common definitions of measurable and simple functions.
\begin{manualdefn}{A.1}
Let $f: X \to \R$ be a measurable function. The \textit{integral} of $f$ is defined to be 
\[\int f \, d\mu = \sup_{\substack{s_f \in \mathcal{S}(X)\\s_f \leq f}} \int s_f \, d\mu\]
where $\mathcal{S}(X)$ is the set of simple functions on $X$.
\end{manualdefn}
\begin{manualtheorem}{A.2}
The integral, as defined above, is a linear operator.
\end{manualtheorem}
\begin{proof}
First, we show linearity for simple functions. That is, let $f$ and $g$ be functions from $X$ to $\R$ that can be written in the form
\[f = \sum_{i=1}^n c_i\bm 1_{S_i} \hspace{1cm} g = \sum_{i=1}^m d_i \bm 1_{T_i}\]
where $m, n < \infty$ are arbitrary, $S_i, T_i \in \Sigma$, and $c_i, d_i \in \R \setminus \{0\}$. Then by definition we have
\[\int kf \, d\mu = \int \sum_{i=1}^n kc_i \bm 1_{S_i} = \sum_{i=1}^n kc_i \mu(S_i) = k\sum_{i=1}^n c_i \mu(S_i) = k \int f \, d\mu\]
And we find
\begin{align*}
\int f + g \, d\mu &= \int \sum_{i=1}^n c_i \bm 1_{S_i} + \sum_{i=1}^m d_i \bm 1_{T_i} \, d\mu \\ &= \sum_{i=1}^n c_i \mu(S_i) + \sum_{i=1}^m d_i \mu(T_i) \\ &= \int \sum_{i=1}^n c_i \bm 1_{S_i} \, d\mu + \int \sum_{i=1}^m d_i \bm 1_{T_i} \, d\mu \\ &= \int f \, d\mu + \int g \, d\mu
\end{align*}
Now we may turn to measurable functions in general. We take as a given that
\[\int cf \, d\mu = c\int f \, d\mu\]
where $\mu$ is finitely additive and $c \in \R$. The proof of this proceeds identically to the proof for countably additive measures. Now let $f$ and $g$ be \textit{measurable non-negative} functions from $X$ to $\R$, and recall the definition of the integral of a non-negative function $f$ as the supremum of the integrals of all nonnegative simple functions bounded by $f$. First observe that for any non-negative simple functions $s_f$ and $s_g$ bounded by $f$ and $g$ respectively, we have
\[s_f \leq f, s_g \leq g \implies s_f + s_g \leq f + g\]
Since the sum of two simple functions is simple, we have
\[\int f + g \, d\mu \geq \int s_f + s_g \, d\mu = \int s_f \, d\mu + \int s_g \, d\mu\]
for all functions $s_f$ and $s_g$, whence
\[\int f + g \, d\mu \geq \sup \left( \int s_f \, d\mu + \int s_g \, d\mu \right) = \int f \, d\mu + \int g \, d\mu\]
The opposite inequality, which is now sufficient to prove additivity of the integral, does not follow as simply, since it is not as easy to split a simple function apart as it is to put two together. Consider a simple function $s_h \leq f + g$. We write
\[s_h = \sum_{i=1}^n c_i \bm 1_{S_i}\]
where, as above, $c_i \in \R \setminus \{0\}$ and $S_i \in \Sigma$. We further assert, for simplicity, that $S_i \cap S_j = \emptyset$ for $i \neq j$.
\begin{manualtheorem}{A.3}
Let $f$ and $g$ be non-negative measurable functions that vanish outside a measurable set $S \in \Sigma$ such that $f + g \geq \bm 1_S$. Then
\[\int f \, d\mu + \int g \, d\mu \geq \mu(S)\]
\end{manualtheorem}
\begin{proof}
Without loss of generality, suppose that $f, g \leq 1$. There is no loss of generality in this supposition because
\[\min(1, f) \leq f \implies \int \min(1, f) \, d\mu \leq \int f \, d\,\mu\]
and
\[f + g \geq \bm 1_S\implies \min(1, f) + \min(1,g) \geq \bm 1_S\]

Now we can place lower bounds on the integrals of $f$ and $g$ by creating sequences $(f^*_n)$ and $(g^*_n)$ of simple functions that are bounded by $f$ and $g$, respectively. In particular, define
\[f^*_n = \sum_{i=0}^{2^n} \bm 1_{f^{-1}((i2^{-n}, (i+1)2^{-n}])} \cdot i2^{-n}\]
which is certainly bounded above by $f$. We define $g^*_n$ similarly based on $g$. We can see that
\[f - f^*_n \leq 2^{-n} \bm 1_S \hspace{1cm} g - g^*_n \leq 2^{-n} \bm 1_S\]
This leads to
\[(f+g) - (f^*_n + g^*_n) \leq 2^{-n+1}\bm 1_S \implies f^*_n + g^*_n \geq (f+g) - 2^{-n+1}\bm 1_S \geq (1-2^{-n+1}) \bm 1_S\]
Now we find
\[\int f \, d\mu + \int g \, d\mu \geq \int f^*_n \, d\mu + \int g^*_n \, d\mu \geq (1-2^{-n+1}) \mu(S)\]
and taking $n \to \infty$ yields
\[\int f \, d\mu + \int g \, d\mu \geq \mu(S)\]
as desired.
\end{proof}
Now we turn back to the main proof. Consider the functions
\[f' = \frac1{c_i} f \cdot \bm 1_{S_i} \hspace{1cm} g' = \frac1{c_i} g \cdot \bm 1_{S_i}\]
we have that $f' + g' \geq \bm 1_{S_i}$, and we need to show that $\int f' \, d\mu + \int g' \, d\mu \geq \mu(S_i)$. This follows directly from the above lemma. This allows us to conclude
\[\int_{S_i} f \, d\mu + \int_{S_i} g \, d\mu \geq c_i \mu(S_i)\]
Since we have assumed the $S_i$ to be disjoint, we can now sum over $i$ to yield that
\[\int f \, d\mu + \int g \, d\mu \geq \int s_h \, d\mu\]
Since this holds for all simple $h \leq f + g$, we find
\[\int f \, d\mu + \int g \, d\mu \geq \sup \left( \int s_h \, d\mu \right) = \int f + g \, d\mu\]
This completes the proof of additivity, and hence the integral is a linear operator.
\end{proof}

\section*{Acknowledgements}
The authors would like to especially thank Qiran Dong, Grace Harper, and Dan Hofman, who were involved with much of the paper's early development. We would also like to thank Dr. William Gasarch and Dr. Clyde Kruskal of the University of Maryland at College Park, as well as the others involved in coordinating the REU program of which this research was originally a part.

Thanks also to Dr. Thomas Goldstein, Dr. Alan Weiss, Yuval Widgerson, and Dr. Wiseley Wong for their comments and suggestions throughout the development of this research.

%\nocite{*}
%\printbibliography
\bibliographystyle{plain}\bibliography{refs}

\end{document}